\documentclass[12pt]{article}
\usepackage[utf8]{inputenc}
\usepackage[T1]{fontenc}
\usepackage{graphicx}
\usepackage[english]{babel}
\usepackage{authblk}
\usepackage{amsmath, amssymb, amsfonts}
\usepackage{geometry}
\usepackage{hyperref}

\begin{document}

\title{Suppression of gel formation by spontaneous fragmentation}
\author[1,2]{Matveev Sergey \thanks{\texttt{matseralex@cs.msu.ru}}}
\affil[1]{Lomonosov Moscow State University, Moscow, Russia}
\affil[2]{Marchuk Institute of Numerical Mathematics RAS, Moscow, Russia}
\date{}
\maketitle

\begin{abstract}
We investigate the discrete Smoluchowski equation with a multiplicative kernel and spontaneous fragmentation at a constant rate $\lambda$. For $\lambda \ge 4$, fragmentation suppresses gelation, and we derive exact expressions for the moments and a recurrence relation for the steady-state distribution. For $0 < \lambda < 4$, we show that post-gelation kinetics depends critically on the closure scheme: the Flory approach (fixed mass) leads to a closure-induced plateau, while the self-consistent Stockmayer approach describes the continuous depletion of the sol. Numerical simulations using the fast Fourier transform-based adaptive Runge-Kutta method fully support the analysis.
\end{abstract}

\section*{Introduction}

The classical Smoluchowski equations provide a fundamental mathematical framework for describing the particle aggregation phenomena \cite{smol1917, galkin2001, tanaka1994}. For the irreversible coagulation process in the spatially homogeneous conditions they describe evolution of the concentrations $n_k(t)$ of particles of sizes $k$ per unit volume and read as an infinite system of differential equations
\begin{align}
\dfrac{d n_k}{dt} 
= \dfrac{1}{2} \sum\limits_{i+j=k} K_{i,j}, n_i n_j 
- n_k \sum\limits_{j=1}^{\infty} K_{i,k} n_j  , \quad k =1, 2 ,3 , \ldots \infty ,
\label{eq:master_Smol}
\end{align}
where $K_{i,j} = K_{j,i} \geqslant 0$ correspond to the coalescence rates (also known as aggregation kernels) between the particles of sizes $i$ and $j$ that form the larger and larger aggregates due to their pairwise collisions $[i]+[j] \rightarrow [i+j]$. One may study such a system of differential equations either analytically or numerically as the Cauchy problem by setting the initial concentrations $n_k(t=0) \geqslant 0$.

A particularly important and exactly solvable case for these equations is the kernel $K_{ij}=i \cdot j$ with monodisperse initial conditions, for which an exact closure of the moment equations is possible \cite{Lushnikov04, lushnikov2006gelation}. In the absence of fragmentation, this kernel leads to gelation, meaning the formation of a macroscopic cluster at a finite time $t_c=1$.

Let us consider the discrete Smoluchowski equation with the product (also known as multiplicative) kernel $K(i,j)=i \cdot j$ and an additional term describing the spontaneous breakup of each particle into monomers at a constant rate $\lambda$:
\begin{align}
\left\lbrace \begin{matrix}
& \dfrac{d n_k}{dt} 
& = & \dfrac{1}{2} \sum\limits_{i+j=k} i\,j\, n_i n_j 
- k\, n_k \sum\limits_{j=1}^{\infty} j n_j  
- \lambda\, n_k , \quad k = 2 ,3 , \ldots \infty \\
& \dfrac{d n_1}{dt} 
& = & - n_1 \sum\limits_{j=1}^{\infty} jn_j  + \lambda \sum\limits_{j=2}^{\infty} j n_j 
\end{matrix} \right.
\label{eq:master_orig}
\end{align}
The aggregation-fragmentation kinetics with the constant aggregation kernel and spontaneous fragmentation has been studied analytically in \cite{srivastava1982simple}. The pure fragmentation equation is addressed in \cite{ziff1985kinetics}. Further analysis with more general rates including demonstrations of the temporal oscillations can be found e.g. in  \cite{pego2020temporal, fortin2023stability, matveev2026generating}. All of these studies correspond to the mass-conserving problems. 

Let $M_1 = \sum\limits_{j=1}^\infty j\, n_j$ denote the total mass of the system, and fix the monodisperse initial conditions $n_k(0)=\delta_{k,1}$, where $\delta_{k,1}$ is the Kronecker symbol. In this case, taking into account the mass conservation, we transform the equations into the following form
\begin{equation}
\left\lbrace \begin{matrix}
& \dfrac{d n_k}{dt} 
& = & \dfrac{1}{2} \sum\limits_{i+j=k} i\,j\, n_i n_j 
- k\, n_k M_1
- \lambda\, n_k , \\
& \dfrac{d n_1}{dt} 
& = & -n_1 M_1 +  \lambda M_1 - \lambda n_1
\end{matrix} \right. .
\label{eq:master_orig_shifted}
\end{equation}
Obviously, for the monodisperse initial conditions $M_1(t) = 1$. Following the mass conservation law \cite{ziff1983kinetics, Flory41}, we substitute this value directly into Eq.~\eqref{eq:master_orig_shifted} and obtain the final form of the target equations
\begin{equation}
\left\lbrace \begin{matrix}
& \dfrac{d n_k}{dt} 
& = & \dfrac{1}{2} \sum\limits_{i+j=k} i\,j\, n_i n_j 
- k\, n_k 
- \lambda\, n_k , \\
& \dfrac{d n_1}{dt} 
& = & \lambda  - n_1 (\lambda + 1) 
\end{matrix} \right. .
\label{eq:master}
\end{equation}

In this work, we investigate the effect of spontaneous fragmentation on the gelation transition with product kernel. We consider the simplest fragmentation mechanism, in which each particle independently breaks down into monomers at a constant rate $\lambda$. For this model, we derive the closed-form analytical expressions for the moments, find a semi-analytical stationary distribution via a generating function and a recurrence relation, and analyze the post-gelation dynamics. 

Similar qualitative observations that spontaneous fragmentation may suppress gelation can also be found in the mathematically oriented literature. In particular, in \cite{laurencot2015absence}
the absence of gelation is proved and  self-similar behavior is established for a continuous coagulation-fragmentation
equation with a multiplicative coagulation kernel and a singular fragmentation kernel. The steady state and stochastic
cyclicity are studied numerically in \cite{fagan2024robustness} for a family of generalized coalescence--fragmentation models with more complicated debris distributions.  

Recent results on the continuous properties of the gelation transition are available for the aggregation processes in the context of random graph growth and explosive percolation in random networks \cite{riordan2011explosive,  riordan2025phase} for the different models including the Erd{\H{o}}s--R{\'e}nyi model \cite{ER60} and Achlioptas process \cite{achlioptas2009explosive}. Analytical investigation of the shattering transition (that is dual to the gelation) recently been addressed in \cite{krapivsky2026aggregation} but for the continuous system with total shatteting of the particles into the dust (that can be also considered as evaporation). 
 
However, to the best of our knowledge, the exact critical fragmentation rate for suppression of gel formation has not been reported previously for the discrete Smoluchowski equation with a multiplicative kernel and constant monomer fragmentation. In this work, we find such an exact critical value $\lambda = 4$ separating the gelling and gel-free regimes, and derive
a closed-form expression for the gelation time $t_c(\lambda)$ for
$0 < \lambda < 4$. 

For the non-gelling regime $\lambda \ge 4$, we establish a recurrence relation that yields the steady-state particle-size distribution with machine precision. Furthermore, for the gelling regime $0 < \lambda < 4$, we provide a detailed comparison of the two common kinetic closures: the Flory approach (fixed total mass) that is given with Eq.~\eqref{eq:master} and the Stockmayer approach (self-consistent finite mass). 

We demonstrate that the Flory closure leads to a stationary state of the closed Flory system that cannot describe the post-gel mass loss, whereas the self-consistent Stockmayer approach correctly captures the continuous absorption of the sol by the gel, reproducing the classical $M_1(t) \sim 1/t$ asymptotics for $\lambda=0$ \cite{ziff1983kinetics}.

\section{Mass conservation and closure of equations}

\subsection{Conditions for gelation}

For standard coagulation kernels (with homogeneity exponent not exceeding one), the formal mass conservation is easy to justify (see, e.g., the book \cite{Krapivsky} and the generalization for the continuous Smoluchowski equation in \cite{galkin2001}). It is sufficient to multiply equations \eqref{eq:master_orig_shifted} by $k$ and sum to obtain $dM_1/dt=0$, whence $M_1(t)\equiv 1$. An important condition for the correctness of such a transformation is that the conditions of the Arzel\`a--Ascoli theorem hold, allowing term-by-term differentiation of the corresponding series. These conditions are satisfied if the homogeneity exponent of the kernel does not exceed one.

Next, we derive an equation for the second moment $M_2 = \sum\limits_{k} k^2 n_k$. Multiplying \eqref{eq:master_orig_shifted} formally by $k^2$ and summing, we obtain
\begin{align}
\frac{d M_2}{dt} 
&= \sum_k k^2 \frac{d n_k}{dt} \nonumber \\
&= \frac12 \sum_{i,j} (i+j)^2 i j n_i n_j - M_1 \sum_k k^3 n_k - \lambda M_2 + \lambda M_1 \nonumber \\
&= (M_3 M_1 + M_2^2) - M_1 M_3 - \lambda M_2 + \lambda M_1,
\end{align}
where \(M_3 = \sum_k k^3 n_k\). The terms containing $M_3$ cancel, and for $M_1=1$ we obtain
\begin{equation}
\frac{d M_2}{dt} = M_2^2 - \lambda M_2 + \lambda, \qquad M_2(0)=1. \label{eq:Riccati}
\end{equation}
This is a Riccati equation with constant coefficients. The case without fragmentation $\lambda = 0$ is well studied \cite{Krapivsky, ziff1983kinetics}. Below we perform a short analysis of the critical time $t_c(\lambda)$ after which gelation begins. At the gelation point the second moment becomes infinite; after passing through this point the mass $M_1$ is not conserved \cite{Flory41,Stockmayer43, Lushnikov04}.

The gelation time is determined from \eqref{eq:Riccati} by direct integration, i.e.
$$
t_c(\lambda) = \int_{1}^{\infty} \frac{dM_2}{M_2^2 - \lambda M_2 + \lambda}.
$$
An analysis of the discriminant of the denominator \(\Delta = \lambda(\lambda-4)\) leads to the simple conclusions:
\begin{itemize}
\item For $\lambda\ge4$, no gel is formed and $M_2(t)$ is bounded;
\item If $0<\lambda<4$ then the critical time $t_c(\lambda)$ is finite and is calculated by the relation
\begin{equation}
t_c(\lambda) = \frac{2}{\sqrt{\lambda(4-\lambda)}} \left[ \frac{\pi}{2} - \arctan\!\left(\frac{2-\lambda}{\sqrt{\lambda(4-\lambda)}}\right) \right];
\end{equation}
\item Separately, the special case $\lambda=0$ is well-known and $t_c=1$.
\end{itemize}
For small $\lambda$, one can obtain the asymptotic expansion for the critical time $t_c(\lambda)=1+\frac{\lambda}{6}+O(\lambda^2)$ by analyzing the corresponding Taylor series. 

\begin{figure}[ht]
\centering
\includegraphics[width=0.49\textwidth]{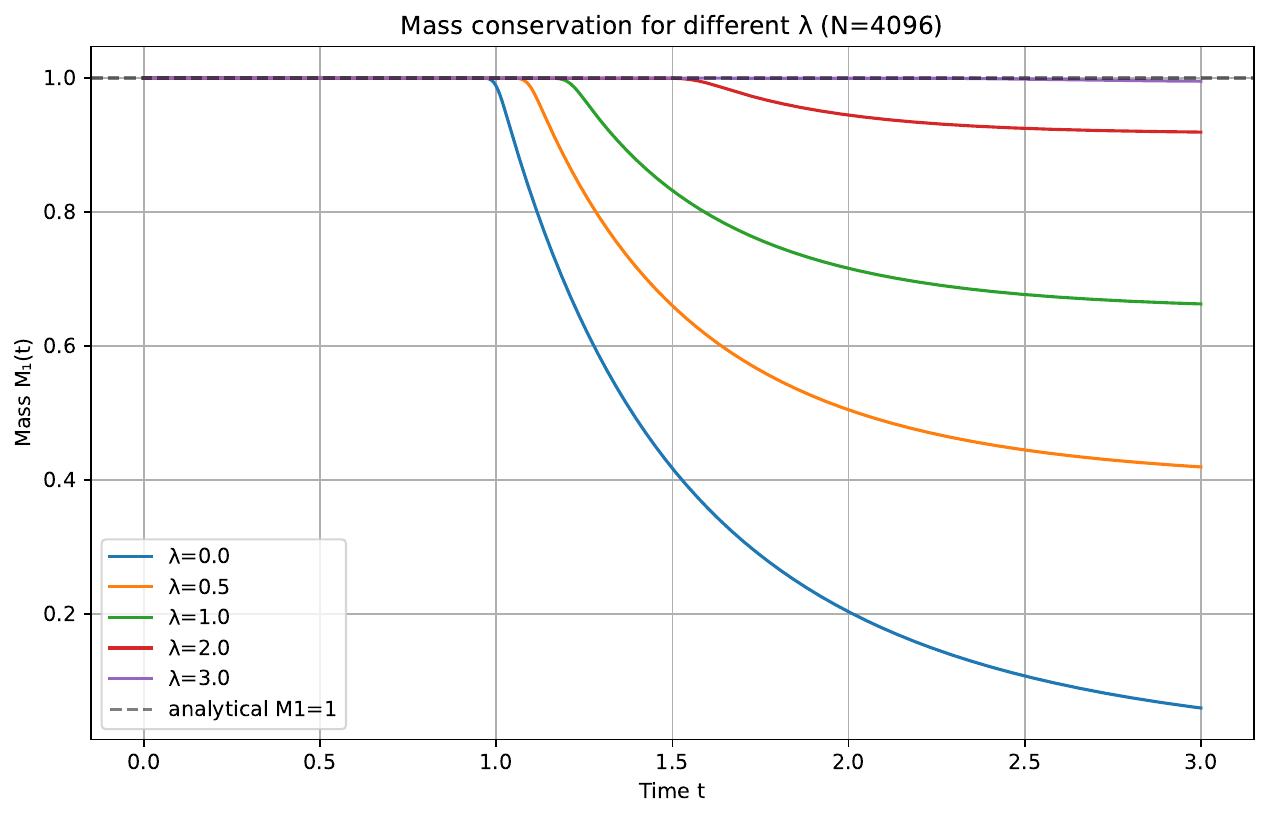}
\includegraphics[width=0.49\textwidth]{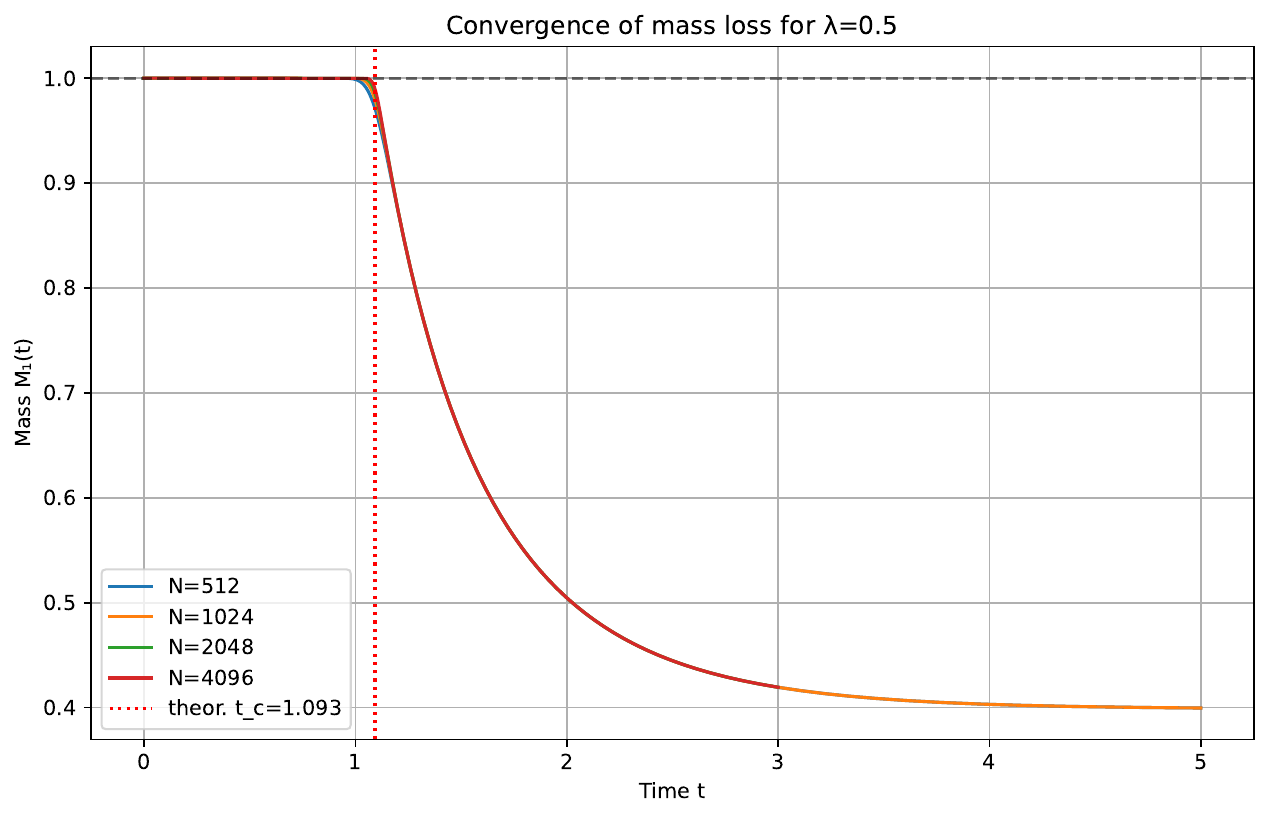}
\caption{Left: dynamics of $M_1(t)$ for different $\lambda$ but with a fixed number of equations $N=4096$. For $\lambda=0$, mass is conserved up to the gelation point and then decreases exponentially. As $\lambda$ increases, gelation occurs later. Right: convergence of the numerical total mass $M_1(t)$ for $\lambda=0.5$ as the number of equations $N$ increases. It is seen that the point at which mass begins to decrease tends to the theoretical $t_c(0.5)$.}
\label{fig:Mass_decrease}
\end{figure}

\subsection{Monomer dynamics for $\lambda \geq 4$}
For $\lambda \geq 4$, gelation does not occur and mass is conserved. Therefore, for our equations with monodisperse initial conditions, for the dynamics of the monomer concentration from \eqref{eq:master} we have the elementary Cauchy problem
$$
\frac{d n_1}{dt} = - (1+\lambda) n_1 + \lambda, \qquad n_1(0)=1,
$$
with the analytical solution
\begin{equation}
n_1(t) = \frac{\lambda}{1+\lambda} + \frac{1}{1+\lambda} e^{-(1+\lambda)t}. \label{eq:n1_sol}
\end{equation}
As $t\to\infty$, the steady-state value is $n_1(\infty)=\frac{\lambda}{1+\lambda} > 0$.

From this it is easy to see that sufficiently intense spontaneous fragmentation prevents the formation of gel that would absorb the finite-size particles. It is interesting that in the case of a system with a product kernel and a constant monomer source but without fragmentation, the total mass decreases as $M_1(t) \simeq 1/t$, i.e., the gel still absorbs the entire mass of the dispersed part of the system \cite{lushnikov2015source, krapivsky2024gelation}.

\begin{figure}[ht]
\centering
\includegraphics[width=0.48\textwidth]{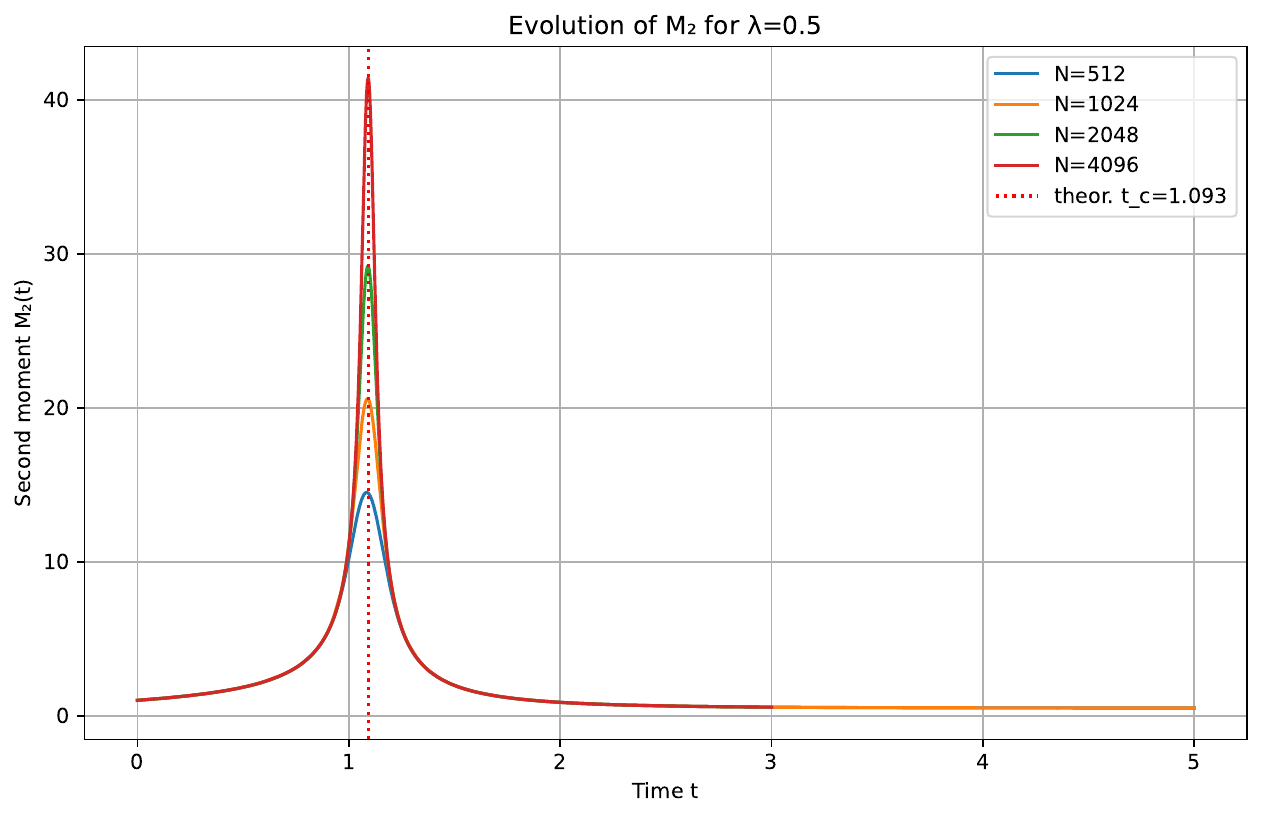}
\includegraphics[width=0.48\textwidth]{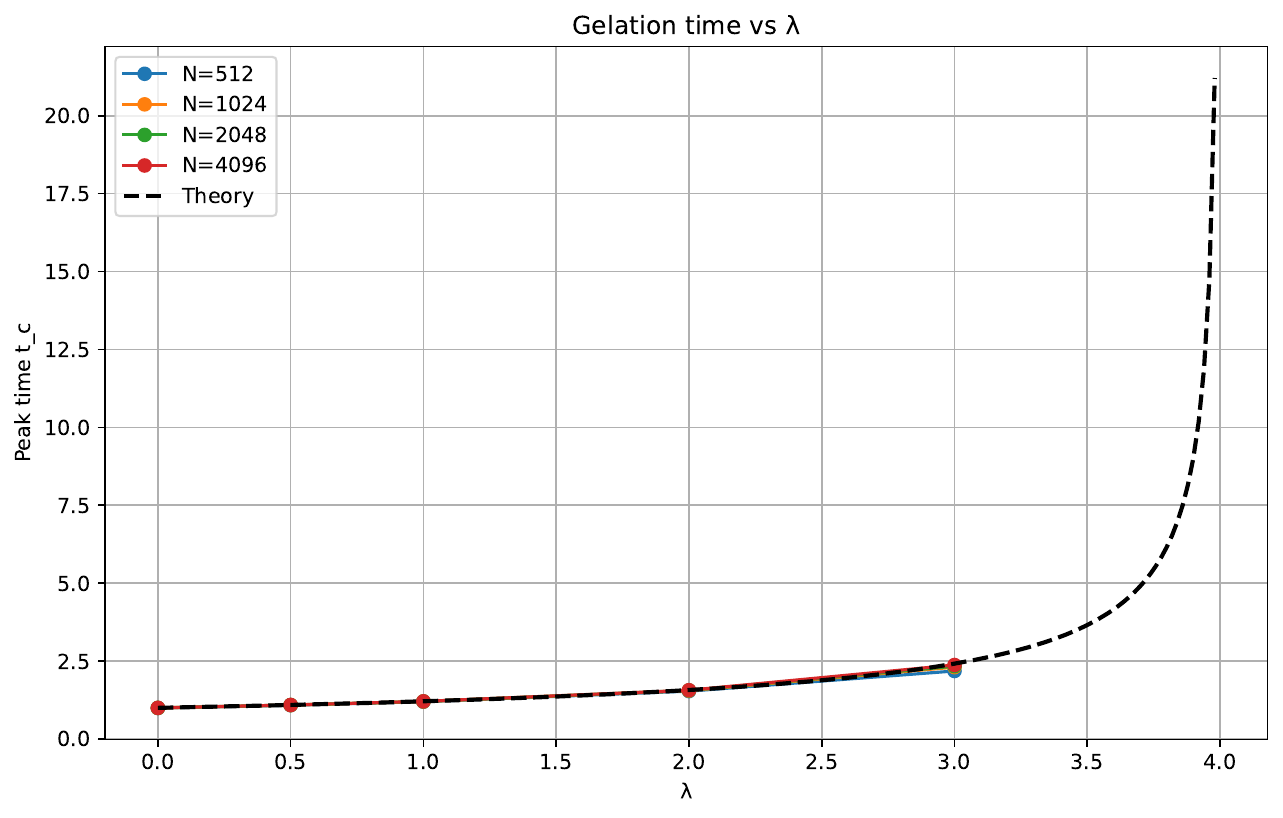}
\caption{Numerical results for the second moment and theoretical estimates. Left: evolution of $M_2(t)$ for $\lambda=0.5$ for different $N$. As $N$ increases, the peak becomes higher and shifts closer toward the theoretical gelation time $t_c(0.5)$. Right: dependence of the peak time $t_c(\lambda)$ on $\lambda$ for different $N$. The dashed line denotes the analytical curve that can be hardly reproduced via numerical integration for $\lambda > 3$.}
\label{fig:M2_dynamics}
\end{figure}

\section{Steady state and generating function}

For $\lambda \ge 4$, the system reaches a steady state and no gel is formed.
For $0 < \lambda < 4$, on the other hand, gelation occurs at a finite time
$t_c(\lambda)$; after gelation the gel continues to absorb finite particles, so
no stationary state with $M_1=1$ exists in this range. This is already visible from the Riccati equation
\eqref{eq:Riccati}, whose fixed points
$M_2 = (\lambda \pm \sqrt{\lambda^2 - 4\lambda})/2$ are real only for
$\lambda \ge 4$.

We therefore analyze the steady state for
$\lambda \ge 4$ only and introduce the generating function for the sequence $n_k$
$$
G(z) = \sum_{k=1}^\infty n_k z^k, \qquad |z|\le 1.
$$
Recall that in the steady state $ d n_k /dt=0$. Multiply \eqref{eq:master} by $z^k$ and sum over $k\ge1$
$$
0 = \frac12 \sum_{k} z^k \sum_{i+j=k} i j n_i n_j - \sum_k k n_k z^k - \lambda \sum_k n_k z^k + \lambda \sum_k \delta_{k,1} z^k.
$$
Using the properties of generating functions
$$
\sum_{k} z^k \sum_{i+j=k} i j n_i n_j = \left( \sum_i i n_i z^i \right)^2 = (z G'(z))^2,
$$
$$
\sum_k k n_k z^k = z G'(z), \qquad \sum_k n_k z^k = G(z), \qquad \sum_k \delta_{k,1} z^k = z,
$$
we obtain the steady-state equation
\begin{equation}
\frac{1}{2} (z G')^2 - z G' - \lambda G + \lambda z = 0. 
\label{eq:G_z_for_recurrent}
\end{equation}
This is a first-order nonlinear differential equation. We rewrite it in $z G'$ as a quadratic equation
$$
(z G')^2 - 2 z G' - 2\lambda G + 2\lambda z = 0.
$$
Solving it, we choose the root containing the minus sign (because $G(0)=0$)
\begin{equation}
z G'(z) = 1 - \sqrt{1 + 2\lambda G(z) - 2\lambda z}. 
\label{eq:G_z_for_moments}
\end{equation}
Further, from eq.~\eqref{eq:G_z_for_moments} one can obtain the moments by differentiating with respect to $z$ and substituting $z=1$. Here $G(1) = \sum n_k = N_\infty$ is the total number of particles in the steady state. Substituting $z=1$ into equation \eqref{eq:G_z_for_moments} and taking into account $G'(1)=M_1=1$, we obtain
\begin{equation}
1 = 1 - \sqrt{1 + 2\lambda N_\infty - 2\lambda} \quad \Rightarrow \quad N_\infty = 1 - \frac{1}{2\lambda}.
\end{equation}
This expression is nonnegative for $\lambda>1/2$. For non-gelling $\lambda\ge4$ it gives $N_\infty > 1/2$ (for example, for $\lambda=5$ one obtains the convenient value $N_\infty = 9/10$).

For the second moment, we differentiate eq.~\eqref{eq:G_z_for_moments} with respect to $z$ and substitute $z=1$. One can show that this gives the following correct root
\begin{equation}
M_2(\infty) = \frac{\lambda - \sqrt{\lambda^2 - 4\lambda}}{2},
\end{equation}
which agrees with the limit of the solution of the Riccati equation as $t\to\infty$.

\begin{figure}[ht]
\centering
\includegraphics[width=0.48\textwidth]{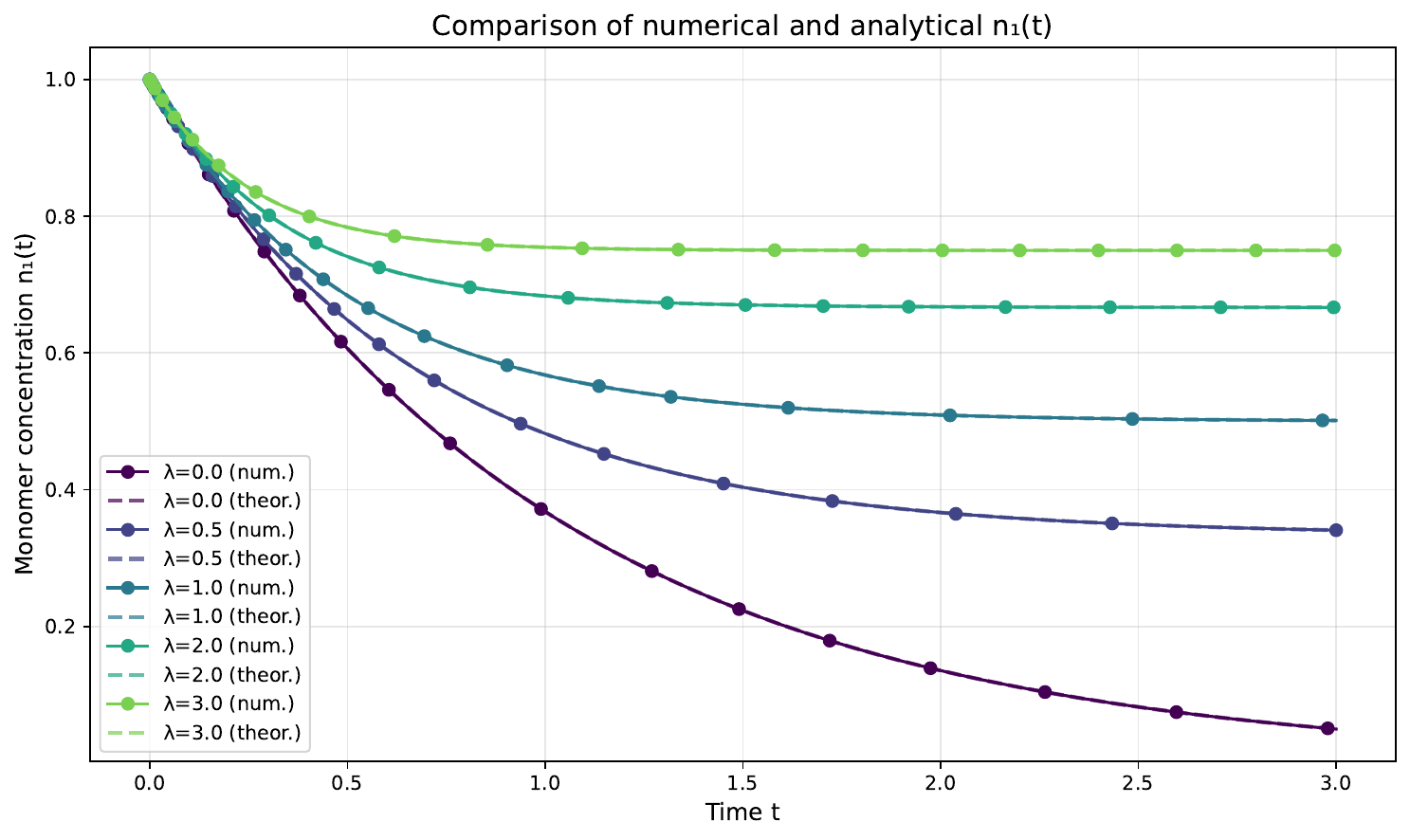}
\includegraphics[width=0.48\textwidth]{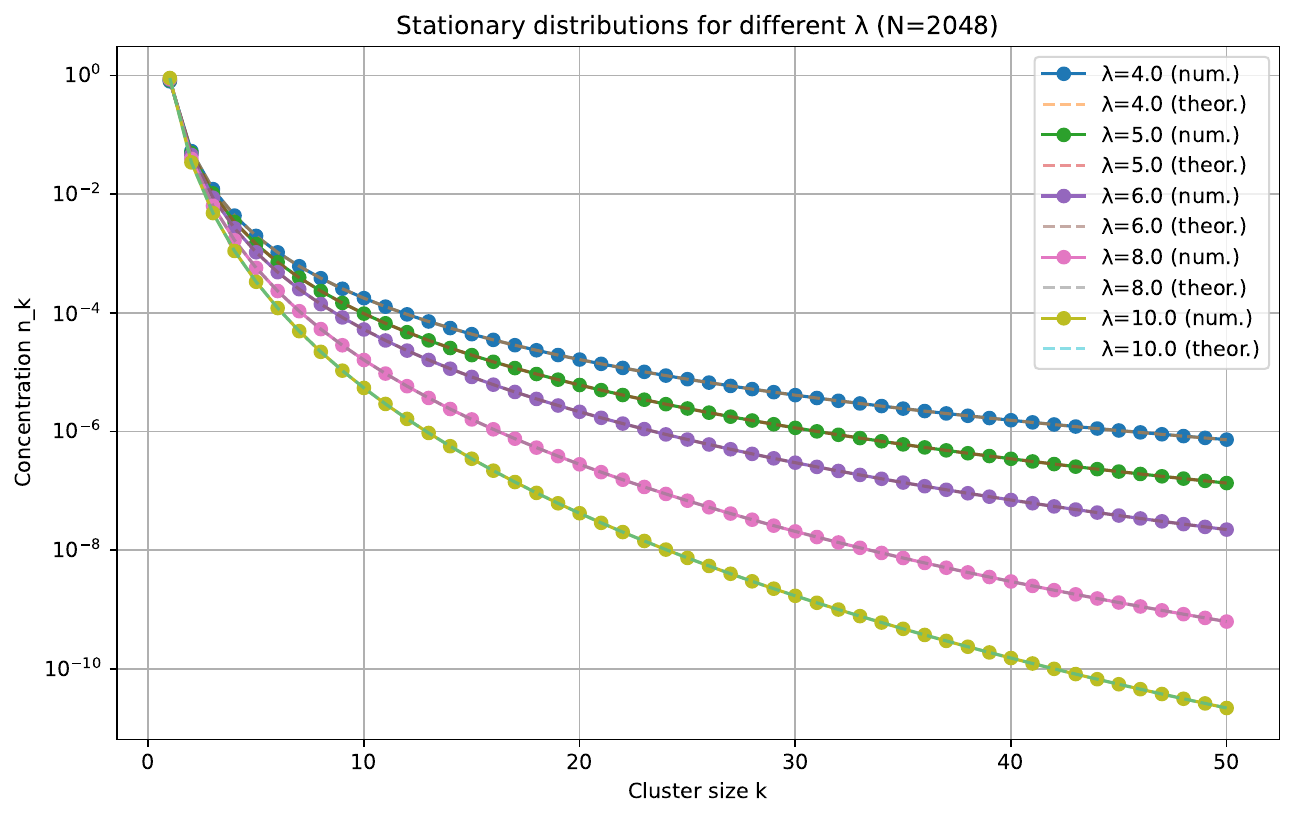}
\caption{Left: comparison of the numerical and semi-analytical solutions for the monomer concentration $n_1(t)$ for various $\lambda$. Right: steady-state distributions $n_k$ for $\lambda \ge 4$ (numerical points and analytical curves from the recurrence relation).}
\label{fig:n1_stationary}
\end{figure}

\section{Recurrence relation for steady-state concentrations}

Instead of solving the equation for the generating function, one can obtain a direct recurrence relation for the concentrations $n_k$. Let us substitute the expansion $G(z) = \sum_{k=1}^\infty n_k z^k$ into the equation for the generating function \eqref{eq:G_z_for_recurrent}. Then $z G'(z) = \sum_{k=1}^\infty k n_k z^k$, and 
$$(z G')^2 = \sum_{k=2}^\infty \left( \sum_{i=1}^{k-1} i(k-i) n_i n_{k-i} \right) z^k.$$
Consequently, the equation \eqref{eq:G_z_for_recurrent} takes the form
$$
\frac{1}{2} \sum_{k=2}^\infty \left( \sum_{i=1}^{k-1} i(k-i) n_i n_{k-i} \right) z^k - \sum_{k=1}^\infty k n_k z^k - \lambda \sum_{k=1}^\infty n_k z^k + \lambda z = 0.
$$
Equating the coefficients of $z^k$ for $k=1$, we obtain the stationary concentration for monomers
    $$-1 \cdot n_1 - \lambda n_1 + \lambda = 0 \quad\Rightarrow\quad n_1 = \frac{\lambda}{1+\lambda},$$
which coincides with \eqref{eq:n1_sol} at $t=\infty$. Further, for $k\ge 2$ we have
$$
\frac{1}{2} \sum_{i=1}^{k-1} i(k-i) n_i n_{k-i} - k n_k - \lambda n_k = 0.
$$
This yields the recurrence formula
\begin{equation} 
n_k = \frac{1}{2(k+\lambda)} \sum_{i=1}^{k-1} i(k-i) n_i n_{k-i}, \qquad k\ge 2.
\label{eq:recurrence}
\end{equation}
This is a convenient result: knowing only the steady-state value for monomers $n_1(t=\infty)$ from equation \eqref{eq:n1_sol} (which depends on $\lambda$), one can successively compute $n_2$, $n_3$ and so on for any $\lambda \ge 4$ with machine precision. Such recurrence forms are rare for aggregation-fragmentation equations and always lead to fruitful theoretical observations \cite{matveev2020oscillating, stadnichuk2015smoluchowski, Krapivsky}. 

For $0<\lambda<4$, the recurrence \eqref{eq:recurrence} with $M_1=1$ yields a normalizable distribution whose mass $\sum_k k n_k$ is strictly less than unity. The difference $1-\sum_k k n_k$ is precisely the mass that, in the original unclosed equation, has already been absorbed by the gel. Thus, the assumption $M_1=1$ that enters \eqref{eq:recurrence} becomes inconsistent with the actual mass of the distribution, which is the mathematical signature of gelation.

For $\lambda \ge 4$, the series for the total number of particles $\sum_k n_k$ converges, and its sum equals $N_\infty = 1 - 1/(2\lambda)$. Thus, the recurrence relation \eqref{eq:recurrence} gives a semi-analytical solution of the steady-state problem, which can be easily implemented numerically to check the results of numerical solvers.

\begin{figure}[ht]
\centering
\includegraphics[width=0.9\textwidth]{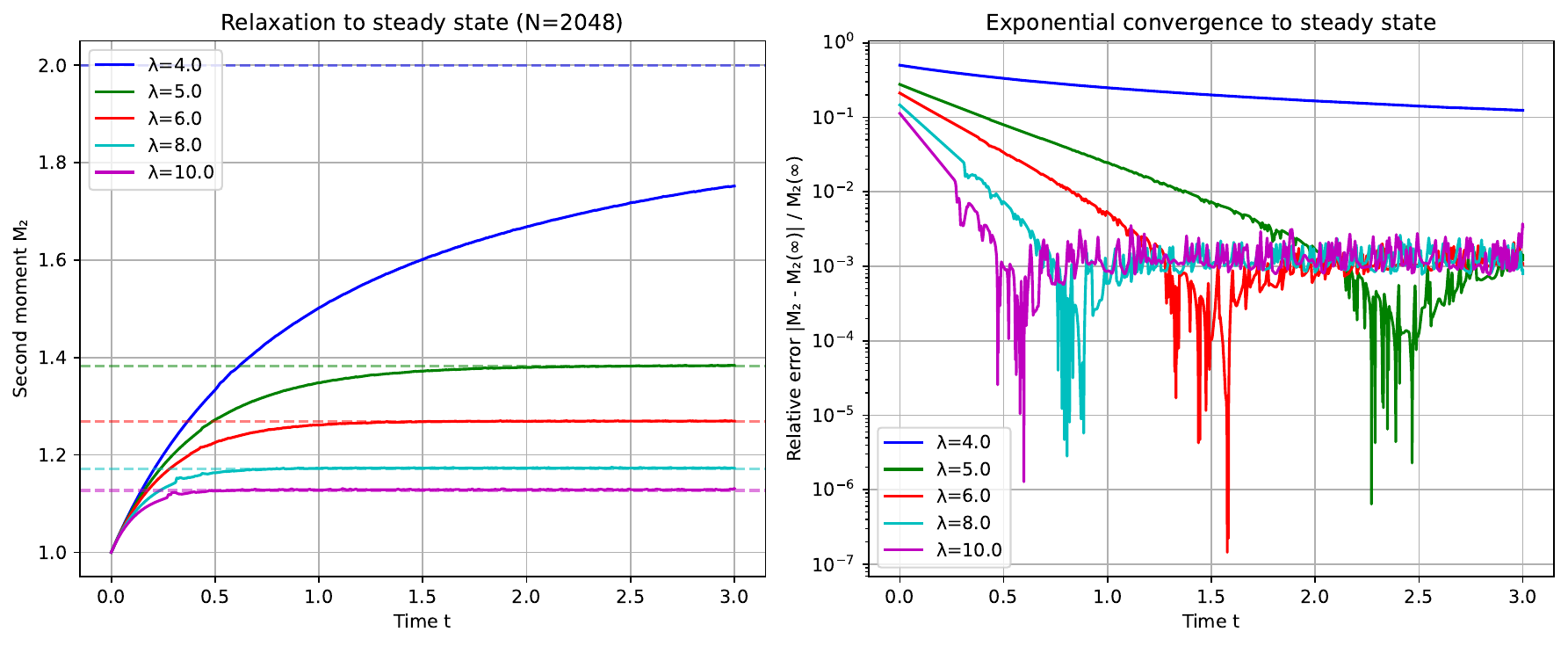}
\caption{Convergence of the second moment to the steady-state value for various $\lambda \ge 4$. Left: $M_2(t)$ (solid lines) and theoretical asymptotics (dashed). Right: relative error on a logarithmic scale, demonstrating exponential decay.}
\label{fig:M2_convergence}
\end{figure}

\section{Comparison with numerical results}

\subsection{General pre-gelation dynamics and steady-state for $\lambda \geq 4$}
\label{sec:validation}

In the numerical solution, we used an algorithm based on the fast Fourier transform to compute the right-hand side in $O(N\log N)$ operations, as well as a second-order Runge-Kutta method with adaptive time-step selection \cite{matveev2015fast, matveev2024adaptive}.

The analytical results can be confirmed numerically: the mass $M_1(t)$ is conserved for $t < t_c$ (see Fig.~\ref{fig:Mass_decrease}), a good quality of approximation is observed for the monomer dynamics $n_1(t)$ from equation \eqref{eq:n1_sol}, and the second moment of the solution $M_2(t)$ exhibits a peak for $\lambda<4$ in the vicinity of $t_c$ (see Fig.~\ref{fig:M2_dynamics}) and steady-state behavior for $\lambda\ge4$.

The quality of the numerical approximation of $t_c(\lambda)$ through the peak value of $M_2(t)$ decreases as $\lambda$ grows and approaches the critical value $\lambda=4$. To perform calculations with good quality, one has to take an increasingly large number of equations; their number grows, which forces one to work with extremely small time steps. In this work, for simplicity of calculations, we restricted ourselves to intervals $t \le 5$, $N \leq 4096$.

The steady-state concentrations computed via recurrence relation \eqref{eq:recurrence} coincide with the numerical ones within the discretization error. At the same time, as $\lambda$ grows, one again has to solve increasingly stiff systems of differential equations. However, the number of equations $N$ required for numerical mass conservation gradually decreases (see Fig. \ref{fig:n1_stationary}).

The steady-state values of $M_2(t)$ for $\lambda \geq 4$ are reached rapidly as the fragmentation intensity grows (see Fig. \ref{fig:M2_convergence}), and the quality of their approximation based on numerical integration by the Runge-Kutta method increases.

\subsection{Post-gelation dynamics for $0 < \lambda < 4$: Flory vs. Stockmayer}
\label{sec:flory_vs_stockmayer}

While for $\lambda \ge 4$ the system reaches a stable stationary state, for $0 < \lambda < 4$ gelation still occurs and no stationary distribution with $M_1=1$ exists. From a numerical perspective, the post-gelation stage is highly sensitive to the chosen truncation and closure scheme \cite{galkin2001, krapivsky2024gelation, ball2011instantaneous, matveev2026generating, lepek2026mini}. Specifically, one can distinguish between two common approaches.

\begin{figure}
\includegraphics[width=0.95\textwidth]{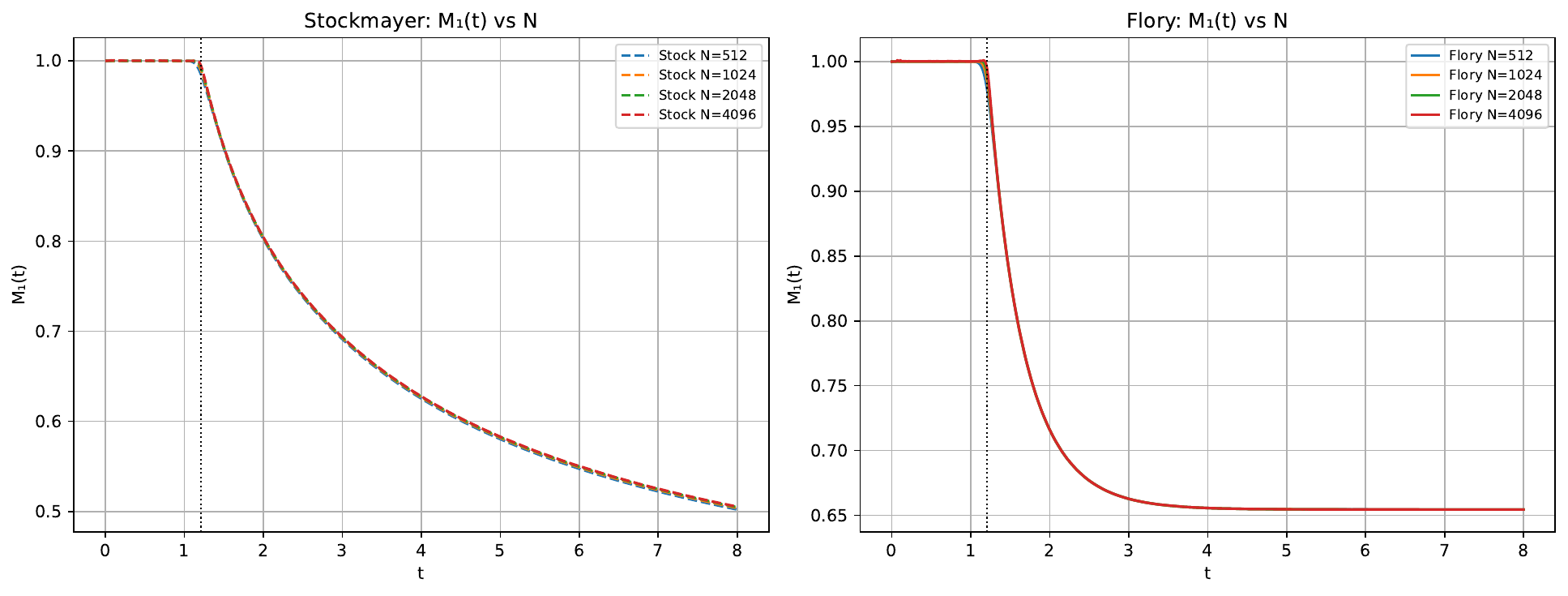}
\caption{Mass dynamics for S-model (left) and F-model (right). In both models the mass decreases after the $t_c$ but in the Flory framework it reaches a stable stationary value, whereas in the Stockmayer closure the gel continues absorption of the finite fraction.}
\label{fig:N_mass_lambda1}
\end{figure}

The first is the Flory closure with direct substitution of $M_1 \equiv 1$, giving Eq.~\eqref{eq:master}. It yields an $N$-stable stationary distribution for any $\lambda > 0$ (see Fig.~\ref{fig:N_mass_lambda1}). For $\lambda \ge 4$ this stationary state coincides with the exact solution of the original Smoluchowski equation $M_1 = 1$. However, for $0 < \lambda < 4$ the stationary state of the F-model has $M_1 < 1$ and satisfies $M_2 = \lambda$. It is a solution of the closed Flory system, not of the original unclosed equation, because the closure artificially uses $M_1 = 1$ in the loss term while the actual mass is $M_1 < 1$. 

The second way to study the dynamics after the critical point $t_c$ is the Stockmayer approach. Here, the mass is tracked self-consistently using the finite sum $M_1^{(N)}(t) = \sum\limits_{k=1}^{N} k n_k(t)$. Thus, the kinetic equations take the form
\begin{equation}
\left\lbrace \begin{matrix}
& \dfrac{d n_k}{dt} = \dfrac{1}{2} \sum\limits_{i+j=k} i\,j\, n_i n_j - k\, n_k M_1^{(N)} - \lambda\, n_k , \\
& \dfrac{d n_1}{dt} = -n_1 M_1^{(N)} + \lambda M_1^{(N)} - \lambda n_1 .
\end{matrix} \right.
\label{eq:stockmayer}
\end{equation}
This approach captures the ongoing absorption of the sol by the gel: as mass leaks beyond the limiting size $k=N$, its truncated value $M_1^{(N)}$ decreases, slowing down the loss term but never allowing the sol fraction to stop decaying. The mass cutoff $k^*(t)$ grows rapidly in this case after the critical point. Hence, even though a finite truncation limit is applied, the dynamics remains stable with respect to growing $N$ (see Fig. \ref{fig:N_mass_lambda1}).

The self-consistent Stockmayer truncation, which uses the actual mass, does not possess a stationary state with $M_1=1$ and instead describes the continuous absorption of the sol by the gel. At $\lambda = 0$, both closures give $M_1 \to 0$, but with different rates: exponential for Flory and $1/t$ for Stockmayer \cite{Lushnikov04, Krapivsky}.

This distinction is crucial. In the input-driven aggregation model without fragmentation \cite{krapivsky2024gelation}, the situation is exactly the opposite: the Flory approach predicts the total depletion of the sol, whereas the Stockmayer approach yields a non-trivial stationary distribution. The presence of spontaneous fragmentation reverses this behavior.

\begin{figure}[ht]
\centering
\includegraphics[width=0.98\textwidth]{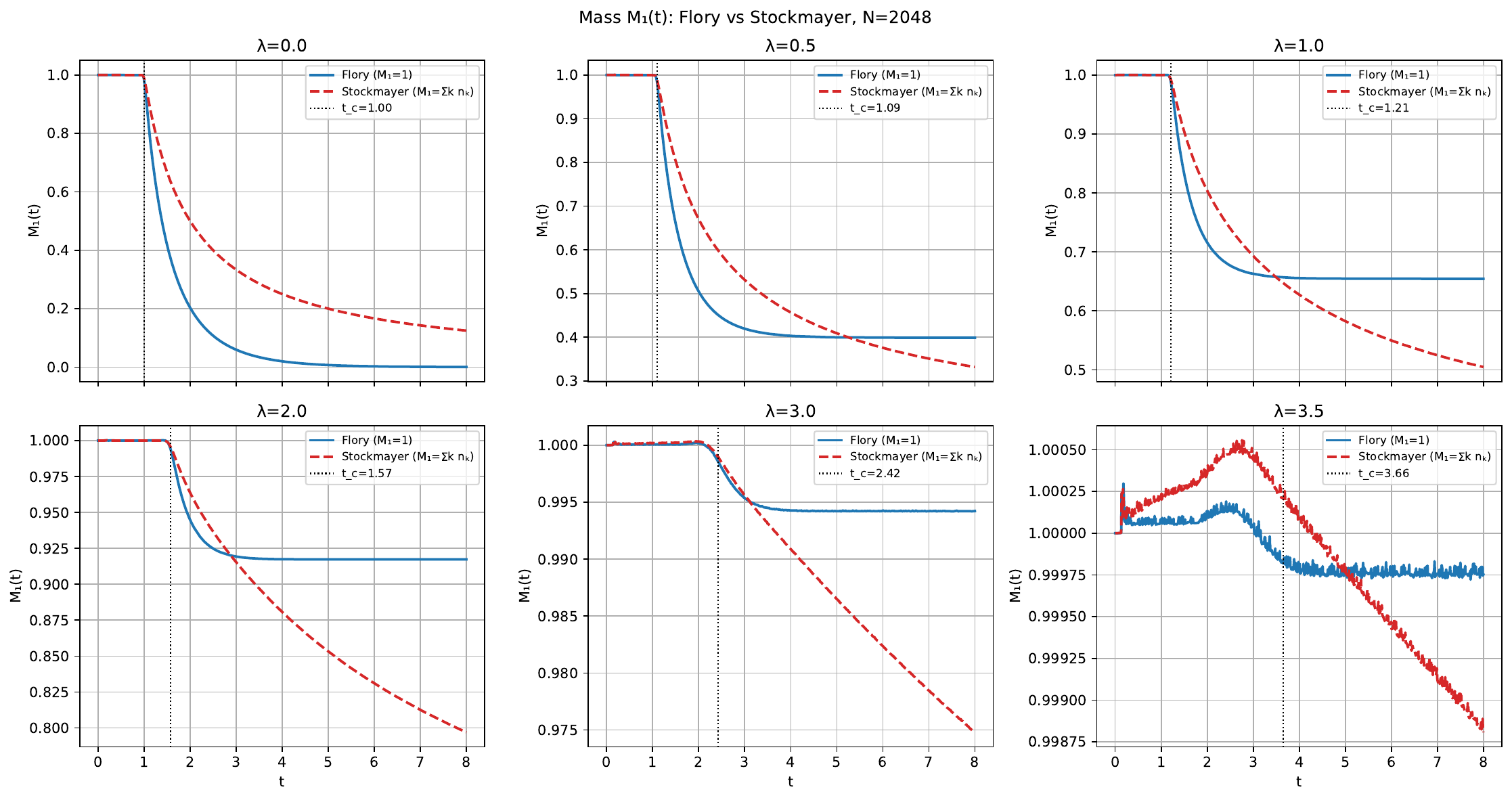}
\caption{Post-gelation dynamics of the sol mass $M_1(t)$ for $N=2048$: comparison of the Flory approach (fixed $M_1=1$, solid blue lines) and the Stockmayer approach (self-consistent $M_1^{(N)}$, dashed red lines) across various fragmentation rates $\lambda \le 3.5$. Vertical dotted lines indicate the theoretical gelation time $t_c(\lambda)$.}
\label{fig:compare_M1}
\end{figure}

In Fig.~\ref{fig:compare_M1} we compare the two approaches across a range of fragmentation rates from $\lambda=0$ up to $\lambda=3.5$. For $\lambda=0.5$, the Flory approach demonstrates a clear plateau $M_1 \approx 0.4$, suggesting a spurious equilibrium between the gel and the sol. Conversely, the Stockmayer approach shows that $M_1^{(N)}(t)$ continues to decrease monotonically. As $\lambda$ increases towards $4$, the plateau in the Flory model persists and its value approaches unity. The Stockmayer curve also decays increasingly slowly, reflecting the divergence of the gelation time $t_c(\lambda)$. The Stockmayer approach with finite $N$ reproduces the power-law tail almost perfectly, whereas the Flory closure leads to a much faster exponential decay (see Fig.~\ref{fig:heavy_tail}).

As model time approaches the critical point $t_c$ for $\lambda \to 4^-$, the system becomes extremely stiff. This is particularly visible in the bottom-right panel of Fig.~\ref{fig:compare_M1} for $\lambda=3.5$, where the adaptive solver struggles to maintain accuracy without highly restrictive time steps, resulting in high-frequency numerical noise and poor quality capture of the gelation point $t_c(\lambda)$. However, despite these numerical artifacts, the qualitative difference between Flory and Stockmayer approaches remains robust: the Flory solution stays on the plateau, while the Stockmayer approach correctly predicts the ongoing decay of the sol.

Finally, we should note a terminological and experimental nuance regarding the terms ``Flory approach'' and ``Stockmayer approach'' \cite{Flory41, Stockmayer43}. In the polymer physics literature, the classical Flory--Stockmayer theory is primarily a statistical framework used to predict the gelation point (the critical conversion degree). Its experimental confirmations typically refer to step-growth polymerization \cite{Adzima2008} and vulcanization \cite{Derosa1998, Kim2012}, where the critical threshold aligns well with theoretical predictions. 

However, concerning the post-gelation asymptotics, the kinetic formalisms diverge. The F-model (frozen mass) and S-model (self-consistent mass) were clearly distinguished for $\lambda=0$ in \cite{Ziff80, ziff1983kinetics}. While the exponential decay of the F-model might be effectively realized in highly constrained polymer networks, the power-law decay of the S-model might be more appropriate for diffusion-limited colloidal aggregation and continuous sol--gel interactions.

\subsection{Critical distribution and the $k^{-5/2}$ tail}
\label{sec:heavy_tail}

\begin{figure}[ht]
\centering
\includegraphics[width=0.9\textwidth]{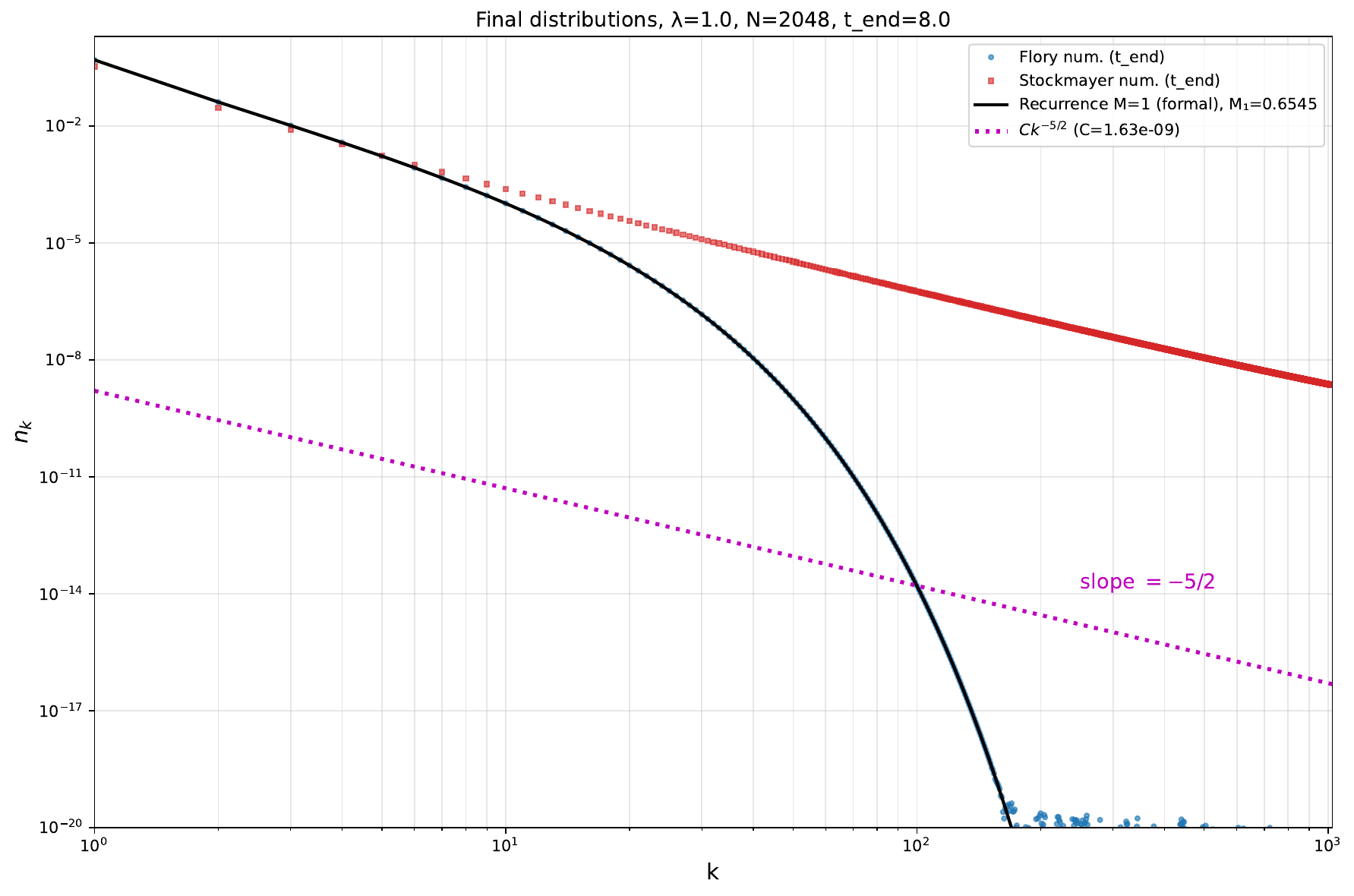}
\caption{Final particle-size distributions for $\lambda=1$ and $N=2048$, at
$t=8$, on a double logarithmic scale. Numerical solution for Flory model corresponds to the blue circles and numerical solution for Stockmayer model is presented with red squares. Solid black line: formal
recurrence \eqref{eq:recurrence} with $M=1$ (the Flory closure).
Dotted violet line: reference power law $C k^{-5/2}$. The numerical Stockmayer distribution follows this power law without exponential cutoff, whereas the Flory solution decays exponentially at large $k$.}
\label{fig:heavy_tail}
\end{figure}

For $0<\lambda<4$ and $t>t_c(\lambda)$, the numerical Stockmayer solution exhibits a power-law tail that is the mark of the incipient gel.
Figure~\ref{fig:heavy_tail} presents the distributions at $t=8$ for $\lambda=1$ and $N=2048$ computed in both closures, together with the formal stationary recurrence \eqref{eq:recurrence} and the reference power law $C k^{-5/2}$.

The formal Flory recurrence with $M_1=1$ (black solid line) gives a distribution whose actual mass equals $M_1 \approx 0.6545<1$. This is inconsistent with the assumption $M_1=1$ underlying the closure: the Flory stationary state does not preserve the mass that in the original unclosed equation has already been absorbed by the gel. For $\lambda=1$ the Flory stationary distribution decays exponentially at large $k$, so its second moment is finite.

In contrast, the numerical Stockmayer distribution (red squares) follows, over two decades in $k$ and without an exponential cutoff up to $k\sim N$, the pure power law
\begin{equation}
n_k \sim C(t)\, k^{-5/2}, \qquad k \to \infty.
\label{eq:heavy_tail}
\end{equation}
The exponent $-5/2$ is the standard Flory--Stockmayer value for the critical (incipient-gel) distribution of the multiplicative kernel \cite{lushnikov2006gelation,  krapivsky2024gelation, Ziff80, ziff1983kinetics}. The mass sum $\sum_k k\,n_k \sim \sum_k k^{-3/2}$ converges, so the sol mass $M_1^{(\mathrm{sol})}$ remains finite, while the second moment $M_2=\sum_k k^2 n_k \sim \sum_k k^{-1/2}$ diverges, which is the standard signature of the gelation threshold.

As $t$ increases, the total sol mass decreases and the prefactor $C(t)$ decreases, consistent with the continuous absorption of the sol by the gel. In the joint limit $N \to \infty$ and $t \to \infty$, the cutoff size moves to infinity and the sol mass vanishes, leaving only the gel.

\section{Conclusion}

An analytical and numerical study of the discrete Smoluchowski equation with a product kernel and the simplest size-independent particle fragmentation into monomers has been carried out. Our analysis includes the derivation of the Riccati equation, the exact solution for $n_1(t)$, and the analysis of gelation and the steady state via the generating function and the recurrence relation.

Exact results for the gelation time $t_c(\lambda)$ have been obtained. For $\lambda \ge 4$, steady-state values of the second moment $M_2(\infty)$, the total number of particles $N_\infty$, and a rule for the recurrence computation of the steady-state particle-size distribution have been obtained. Numerical calculations fully confirm the analytical results and the obtained relation, which allows one to compute the steady-state concentrations $n_k$ with machine precision.

Furthermore, we analyzed the post-gelation dynamics for $0 < \lambda < 4$ and demonstrated that the choice of the closure scheme is crucial. We showed that the commonly used Flory closure ($M_1 \equiv 1$) leads to a stationary state of the closed Flory system that does not describe the post-gel mass loss, whereas the self-consistent Stockmayer approach ($M_1 = \sum k n_k$) captures the continuous absorption of the sol by the gel, ensuring $M_1(t) \to 0$ as $t \to \infty$. 

Thus, we can conclude that (1) intense spontaneous fragmentation indeed
suppresses gelation for $\lambda \ge 4$, so that a genuine stationary
distribution exists and can be computed via the recurrence rule;
(2) for $0<\lambda<4$ gelation still occurs, the gelation time $t_c(\lambda)$
grows with $\lambda$ and diverges as $\lambda \to 4^-$, and after $t_c$ the
gel continues to absorb the sol, so no stationary state with $M_1=1$ is reached. 

\section{Acknowledgements}

This work was supported by the Moscow Center of Fundamental and Applied Mathematics at INM RAS (Agreement with the Ministry of Science and Higher Education of the Russian Federation No.075-15-2025-347). 

This scientific work is dedicated to the memory of Alexey A. Lushnikov, whose contribution to the development of the theory of aggregation-fragmentation kinetics and the sol-gel transition can hardly be overstated.

\bibliographystyle{unsrt}
\bibliography{main}

\end{document}